\documentclass{article}
\usepackage{a4wide}%これで横幅を広げている
\usepackage[title]{appendix}
\usepackage{amsthm}%定理環境
\newtheorem{Definition}{Definition}[section]
\newtheorem{Theorem}[Definition]{Theorem}
\newtheorem{Lemma}[Definition]{Lemma}

\newtheorem{Corollary}[Definition]{Corollary}
\newtheorem*{Proof}{Proof}
\newtheorem*{proof_main}{Proof of Corollary1.3}
\newtheorem*{proof_mainT}{Proof of Theorem1.2}

\usepackage{amsmath}%数式環境 alignなど
\DeclareMathOperator{\low}{Low}
\usepackage{amssymb}%数学記号 mathbなど
\usepackage{comment}

\usepackage{enumerate}%箇条書き
\usepackage{authblk}
\usepackage{hyperref}
\newcommand*{\email}[1]{%
    \normalsize\href{mailto:#1}{#1}\par
    }
    
\title{No eigenvectors embedded in the singular continuous spectrum of Schr\"odinger operators}
\author{Kota Ujino}
\date{}
\affil{Graduate School of Mathematics, Kyushu University,\\ 744 Motooka Nishi-ku, Fukuoka 819-0395, Japan\\ \email{ujino.kouta.699@s.kyushu-u.ac.jp}}

\begin{document}
\maketitle
\begin{abstract}
It is known that the spectrum of Sch\"{o}dinger operators with sparse potentials consists of singular continuous spectrum. 
We give a sufficient condition so that the edge of the singular continuous spectrum is not an eigenvalue and construct examples with singular continuous spectrum which have no eigenvalues and which have a single negative eigenvalue.
\end{abstract}
\section{Introduction and a result} 
We investigate one-dimensional Schr\"{o}dinger operators with sparse potentials. 
It is known that the spectrum of Schr\"{o}dinger operators with sparse potential consists of singular continuous spectrum.
Simon and Spencer $\cite{SimonSpencer}$ show the absence of absolutely continuous spectrum of Schr\"odinger operators with sparse potentials. 
Simon and Stoltz $\cite{Simonsparse}$ also show that $-\frac{d^2}{dx^2}f+Vf=Ef$ has no  $L^2$-solutions  for any $E>0$.
We have the question whether the edge of the singular continuous spectrum is an eigenvalue or not.
We give a sufficient condition for the absence of embedded eigenvalues and give  examples.

\begin{Definition}
A function $V:[0,\infty)\rightarrow \mathbb{R}$ is called a sparse potential, if there exist positive sequences $\{x_n\}_{n=1}^{\infty}, \{\alpha_n\}_{n=1}^{\infty}$ and $\{h_n\}_{n=1}^{\infty}$ such that $x_{n+1}>x_n$ for $n=1,2,3,...$, 
\begin{enumerate}[$($\rm i$)$]
\item
$\begin{displaystyle}\lim_{n\rightarrow \infty}\frac{x_{n+1}-x_n}{\alpha_{n+1}+\alpha_n+1}=\infty\end{displaystyle}$,
\item
$|V(x)|\leq h_n$, if $x \in [x_n-\alpha_n,x_n+\alpha_n]$ for $n=1,2,3,...$,
\item
$V(x)=0$, if $\begin{displaystyle}  x \in \left(\bigcup_{n=1}^{\infty} [x_n-\alpha_n,x_n+\alpha_n] \right)^c\end{displaystyle}$
\end{enumerate}
\end{Definition}
We define $L_n=x_{n+1}-x_n-\alpha_{n+1}-\alpha_n$ for $n\geq1$ and $L_0=x_1-\alpha_1$. 
By $\rm(i)$, $L_n\rightarrow \infty$ as $n\rightarrow \infty$. 
By Strum-Liouville theory $\cite[{\rm Theorem}\: 9.1.]{Strum}$, there exists a unique solution $f\in AC_{loc}([0,\infty))$ of the equation $-\frac{d^2}{dx^2}f+Vf=0$ with $\frac{d}{dx}f \in AC_{loc}([0,\infty))$ and the boundary condition\\ 
$f(0)=\alpha,\frac{d}{dx}f(0)=\beta,\alpha,\beta\in\mathbb{C}$.
We give a sufficient condition of the existence of a non $L^2$-integrable solution.
\begin{Theorem}\label{main}
Let $V$ be the sparse potential and $f$ a weak solution of $-\frac{d^2}{dx^2}f+Vf=0$.
If 
\begin{eqnarray}\label{assumption}
\cfrac{L_n}{4^n}
\left(
\prod_{m=1}^{n}(L_{m-1}^2+2)
\right)^{-1}
\left(
\prod_{m=1}^{n}(2\alpha_m^2+1)
\right)^{-1}
\exp
\left(-\frac{2}{3}\sum_{m=1}^{n}h_m(4\alpha_m^3+3\alpha_m)
\right)
\rightarrow \infty,
\end{eqnarray}
as $n\rightarrow \infty$, then $f\notin L^2([0,\infty))$.
\end{Theorem}

We give an example for one-dimensional Schr\"{o}dinger operators with singular continuous spectrum which has no embedded eigenvalues.
Let $x_n=\exp(n^n)$ for $n=1,2,3,...$, and 
\begin{eqnarray*}
V(x)=
\begin{cases}
e^{n},&\text{if $|x-x_n|\leq\frac{1}{2}$ for $n=1,2,3,...$},\\
0,&\text{otherwise}.
\end{cases}
\end{eqnarray*}
Let $H=-\frac{d^2}{dx^2}+V:L^2([0,\infty))\rightarrow L^2([0,\infty))$ with the domain $D(H)=C^{\infty}_0((0,\infty))$. 
We see that $H$ is regular at zero and in the limit point case at infinity.
This implies that $H$ has self-adjoint extensions $H_\theta$ which can be parametrized by boundary conditions. 
Hence, $H_{\theta}$ is the restriction of $H^*$ to $D_{\theta}=\{f\in D(H^*)\mid {f(0)}\sin\theta-{\frac{d}{dx}f(0)}\cos\theta=0\}$. 
By $\cite{Simonsparse}$, we have  $\sigma_{sc}(H_\theta)=[0,\infty)$, $\sigma_{pp}(H_\theta)\cap(0,\infty)=\emptyset$ and $\sigma_{ac}(H_\theta)=\emptyset$ for all $\theta\in (\frac{\pi}{2},\frac{\pi}{2}]$. 
See Appendix A for the proof.
Theorem $\ref{main}$ implies the next corollary.
This also implies that $H_\theta$ has purely singular continuous spectrum for some $\theta$.
\begin{Corollary}
It follows that 
\begin{enumerate}[$(1)$]
\item
$\sigma_{sc}(H_\theta)=[0,\infty)$, $\sigma_{pp}(H_\theta)\cap[0,\infty)=\emptyset$ and $\sigma_{ac}(H_\theta)=\emptyset$ for all $\theta\in(-\frac{\pi}{2},\frac{\pi}{2}]$,
\item
$H_{\theta}$ has purely singular continuous spectrum for $\theta\in[0,\frac{\pi}{2}]$,
\item
$H_{\theta}$ has a single negative eigenvalue for 
$\theta\in (-\frac{\pi}{2},\arctan(-\frac{1+\sqrt{3}}{2})]$.
\end{enumerate}
\end{Corollary}
\section{Proof of Theorem 1.2.}
We calculate a lower bound of Wronskian matrices. 
\begin{comment}
Let  $
\left|
\left(
\begin{array}{c}
a\\
b
\end{array}
\right)
\right|=\sqrt{a^2+b^2}$ and \\
$
\left(
\begin{array}{c}
a\\
b
\end{array}
\right)\cdot
\left(
\begin{array}{c}
c\\
d
\end{array}
\right)
=ac+bd$ for $a,b,c,d\in \mathbb{R}$. 
\end{comment}
For a $2\times2$-real matrix $M$, let 
\[
\low M=
\inf
\left\{
\left|
M
\begin{pmatrix}
\cos \theta\\
\sin \theta
\end{pmatrix}
\right| \:
\middle|\:
\theta \in [0,2\pi)
\right\}.
\]
We see that $|Mu|\geq \low M |u|$ for $u \in \mathbb{R}^2$ and $\low M=\sqrt{\inf\sigma({}^t\!MM)}$.
Let $V$ be a sparse potential and $f\in AC_{loc}([0,\infty))$ be a weak solution of $-\frac{d^2}{dx^2}f+Vf=0$ with $\frac{d}{dx}f\in AC_{loc}([0,\infty))$ and the boundary condition $f(0)=\cos\theta,\frac{d}{dx}f(0)=\sin\theta$. 
We can represent the weak solution concretely as follows. 
Let $J_0=[0,x_1-\alpha_1]$, $I_n=[x_n-\alpha_n,x_n+\alpha_n]$ and $J_n=[x_n+\alpha_n,x_{n+1}-\alpha_{n+1}]$ for $n=1,2,...$. 
Let $p_n, q_n :J_n \rightarrow \mathbb{R}$ be defined by $p_n(x)=1, q_n(x)=x-x_n-\alpha_n$, and $p_0, q_0 :J_0 \rightarrow \mathbb{R}$ by $p_0(x)=1, q_0(x)=x$. 
Let $\varphi_n,\psi_n:I_n\rightarrow \mathbb{R}$ be the weak solution of $-\frac{d^2}{dx^2}f+Vf=0$ on $I_n$ with the boundary condition $\varphi_n(x_n-\alpha_n)=1,\frac{d}{dx}\varphi_n(x_n-\alpha_n)=0, \psi_n(x_n-\alpha_n)=0$ and $\frac{d}{dx}\psi_n(x_n-\alpha_n)=1$.
We define $\tilde{f}:[0,\infty)\rightarrow \mathbb{R}$ by 
\begin{eqnarray*}
\tilde{f}(x)=c_n^{(1)}p_n(x)+c_n^{(2)}q_n(x) \qquad \text{if $x\in J_n$},\\
\tilde{f}(x)=d_n^{(1)}\varphi_n(x)+d_n^{(2)}\psi_n(x)\qquad \text{if $x\in I_n$},
\end{eqnarray*}
where $c_n^{(1)},c_n^{(2)},d_n^{(1)}$ and $d_n^{(2)}$ are inductively determined by for $n=1,2,...$
\begin{eqnarray}\nonumber
\lim_{x\downarrow 0}
\begin{pmatrix}
p_{0}(x)&q_{0}(x)\\
\frac{d}{dx}p_{0}(x)&\frac{d}{dx}q_{0}(x)
\end{pmatrix}
\begin{pmatrix}
c_0^{(1)}\\
c_0^{(2)}
\end{pmatrix}
&=&
\begin{pmatrix}
\cos \theta \\
\sin \theta
\end{pmatrix},
\\\label{c1}
\lim_{x\uparrow x_n-\alpha_n}
\begin{pmatrix}
p_{n-1}(x)&q_{n-1}(x)\\
\frac{d}{dx}p_{n-1}(x)&\frac{d}{dx}q_{n-1}(x)
\end{pmatrix}
\begin{pmatrix}
c_{n-1}^{(1)}\\
c_{n-1}^{(2)}
\end{pmatrix}
&=&
\lim_{x\downarrow x_n-\alpha_n}
\begin{pmatrix}
\varphi_{n}(x)&\psi_{n}(x)\\
\frac{d}{dx}\varphi_{n}(x)&\frac{d}{dx}\psi_{n}(x)
\end{pmatrix}
\begin{pmatrix}
d_{n}^{(1)}\\
d_{n}^{(2)}
\end{pmatrix},
\\\label{c2}
\lim_{x\uparrow x_n+\alpha_n}
\begin{pmatrix}
\varphi_{n}(x)&\psi_{n}(x)\\
\frac{d}{dx}\varphi_{n}(x)&\frac{d}{dx}\psi_{n}(x)
\end{pmatrix}
\begin{pmatrix}
d_{n}^{(1)}\\
d_{n}^{(2)}
\end{pmatrix}
&=&
\lim_{x\downarrow x_n+\alpha_n}
\begin{pmatrix}
p_{n}(x)&q_{n}(x)\\
\frac{d}{dx}p_{n}(x)&\frac{d}{dx}q_{n}(x)
\end{pmatrix}
\begin{pmatrix}
c_{n}^{(1)}\\
c_{n}^{(2)}
\end{pmatrix}
.
\end{eqnarray}
By the definition, $\tilde{f}$ and $\frac{d}{dx}\tilde{f}$ are continuous, and  
$\left(
\begin{array}{c}
\tilde{f}(0)\\
\frac{d}{dx}\tilde{f}(0)
\end{array}
\right)
=
\left(
\begin{array}{c}
\cos \theta\\
\sin \theta
\end{array}
\right)$. 
The coefficients $c_n^{(1)},c_n^{(2)},d_n^{(1)}$ and $d_n^{(2)}$ are uniquely determined by $\theta$.
We see that 
\[
\int_{0}^{\infty}\tilde{f}(x)\frac{d^2}{dx^2}g(x)dx=\int_{0}^{\infty}\tilde{f}(x)V(x)g(x)dx
\]
for $g \in C_0^{\infty}((0,\infty))$ straightforwardly.
This implies that $\tilde{f}$ satisfies $-\frac{d^2}{dx^2}f+Vf=0$ in the sence of the weak derivative. 
Therefore $f=\tilde{f}$ by the uniqueness of the solution.

\begin{Lemma}
Let $c_n=
\begin{pmatrix}
c_n^{(1)}\\
c_n^{(2)}
\end{pmatrix}
.$
Then we see that
\begin{eqnarray}\label{lower bound of L^2 integral}
\int_{J_n}|f(x)^2|dx
\geq
\frac{1}{4}
\frac{L_n^4}{L_n^3+3L_n}|c_n|^2.
\end{eqnarray} 
\end{Lemma}
\begin{Proof}\rm
We obtain
\begin{eqnarray*}
\int_{J_n}|f(x)|^2dx&=&\int_{J_n}|c_n^{(1)}p_n(x)+c_n^{(2)}q_n(x)|^2dx
\\
&=&
\begin{pmatrix}
c_n^{(1)}&
c_n^{(2)}
\end{pmatrix}
\begin{pmatrix}
L_n&\frac{1}{2}L_n^2\\
\frac{1}{2}L_n^2&\frac{1}{3}L_n^3
\end{pmatrix}
\begin{pmatrix}
c_n^{(1)}\\
c_n^{(2)}
\end{pmatrix}.
\end{eqnarray*} 
The matrix 
$
\begin{pmatrix}
L_n&\frac{1}{2}L_n^2\\
\frac{1}{2}L_n^2&\frac{1}{3}L_n^3
\end{pmatrix}$ can be diagonalized and its eigenvalues $\lambda_{\pm}$ are
\[
\lambda_\pm=\cfrac{1}{2}
\left(
L_n+\frac{1}{3}L_n^3
\pm
\sqrt{(L_n+\frac{1}{3}L_n^3)^2-\frac{1}{3}L_n^4}
\right).
\]
Since $1-\sqrt{1-t}\geq\frac{1}{2}t$ for $0<t<1$, we have 
\begin{eqnarray}\nonumber
\lambda_-=
\cfrac{L_n+\frac{1}{3}L_n^3}{2}
\left(1
-\sqrt{1-\frac{L_n^4}{3}\left(L_n+\frac{1}{3}L_n^3\right)^{-2}}
\right)
%\\ \nonumber
%&\geq&
%\frac{1}{2}\:
%\frac{L_n^4}{3}\left(L_n+\frac{1}{3}L_n^3\right)^{-2}
%\cfrac{L_n+\frac{1}{3}L_n^3}{2}
\geq
\frac{1}{4}
\frac{L_n^4}{L_n^3+3L_n}.
\end{eqnarray}
This implies our assertion.
\qed
\end{Proof}
%It is sufficient to prove that $\lim_{n\rightarrow\infty}\frac{L_n^4}{L_n^3+3L_n}|c_n|^2$..
By $(\ref{c1})$ and $(\ref{c2})$, we obtain $c_n=R_nW_{n-1}c_{n-1}$ for $n\geq1$, where
\begin{eqnarray*}
R_m&=&\lim_{x\uparrow x_m+\alpha_m}
\begin{pmatrix}
\varphi_{m}(x)&\psi_{m}(x)\\
\frac{d}{dx}\varphi_{m}(x)&\frac{d}{dx}\psi_{m}(x)
\end{pmatrix},
\\
W_m&=&
\lim_{x\uparrow x_{m+1}-\alpha_{m+1}}
\begin{pmatrix}
p_{m}(x)&q_{m}(x)\\
\frac{d}{dx}p_{m}(x)&\frac{d}{dx}q_{m}(x)
\end{pmatrix}
=
\begin{pmatrix}
1&L_m\\
0&1
\end{pmatrix}
.
\end{eqnarray*}
We shall estimate $\low R_m$ and $\low W_m$. 
\begin{Lemma}
It follows that 
\begin{eqnarray}\label{lower bound of W_m}
\low W_m\geq \frac{1}{\sqrt{L_m^2+2}}.
\end{eqnarray}
\end{Lemma}
\begin{Proof}\rm
The eigenvalues $\xi_{\pm}$ of ${}^tW_mW_m$ are 
\[
\xi_{\pm}=
\frac{1}{2}\left(
L_m^2+2\pm
\sqrt{(L_m^2+2)^2-4}
\right).
\]
Since $1-\sqrt{1-t}\geq\frac{1}{2}\:t$ for $0<t<1$, we have 
\begin{eqnarray*}
\xi_{-}=\frac{L_m^2+2}{2}\left(
1-
\sqrt{1-4(L_m^2+2)^{-2}}
\right)
\geq
\frac{1}{L_m^2+2}.
\end{eqnarray*}
Since $\low M=\sqrt{\inf\sigma({}^t\!MM)}$, we have our assertion.
\qed
\end{Proof}

%To estimate $\low R_m$, we use variation of constants. 
Let $\tilde{\varphi}_m,\tilde{\psi}_m:I_m\rightarrow \mathbb{R}$ be defined by $\tilde{\varphi}_m(x)=1,\tilde{\psi}_m(x)=x-x_m+\alpha_m$.
We see that $\tilde{\varphi}_m,\tilde{\psi}_m$ satisfy $-\frac{d^2}{dx^2}f=0.$
There exist $u_m^{(j)},v_m^{(j)}\in AC(I_m)$, $j=1,2$ such that 
\begin{eqnarray*}
\begin{pmatrix}
\varphi_m(x)\\
\frac{d}{dx}\varphi_m(x)
\end{pmatrix}
=
u_m^{(1)}(x)
\begin{pmatrix}
\tilde{\varphi}_m(x)\\
\frac{d}{dx}\tilde{\varphi}_m(x)
\end{pmatrix}
+
u_m^{(2)}(x)
\begin{pmatrix}
\tilde{\psi}_m(x)\\
\frac{d}{dx}\tilde{\psi}_m(x)
\end{pmatrix},
\\
\begin{pmatrix}
\psi_m(x)\\
\frac{d}{dx}\psi_m(x)
\end{pmatrix}
=
v_m^{(1)}(x)
\begin{pmatrix}
\tilde{\varphi}_m(x)\\
\frac{d}{dx}\tilde{\varphi}_m(x)
\end{pmatrix}
+
v_m^{(2)}(x)
\begin{pmatrix}
\tilde{\psi}_m(x)\\
\frac{d}{dx}\tilde{\psi}_m(x)
\end{pmatrix}.
\end{eqnarray*}
We see that $u_m^{(1)}(x_m-\alpha_m)=1,u_m^{(2)}(x_m-\alpha_m)=0$, 
$v_m^{(1)}(x_m-\alpha_m)=0,v_m^{(2)}(x_m-\alpha_m)=1$, and 
\begin{eqnarray}\label{derivative of constants}
\begin{pmatrix}
\varphi_{m}&\psi_{m}\\
\frac{d}{dx}\varphi_{m}&\frac{d}{dx}\psi_{m}
\end{pmatrix}
=
\begin{pmatrix}
\tilde{\varphi}_{m}&\tilde{\psi}_{m}\\
\frac{d}{dx}\tilde{\varphi}_{m}&\frac{d}{dx}\tilde{\psi}_{m}
\end{pmatrix}
\begin{pmatrix}
u_{m}^{(1)}&v_{m}^{(1)}\\
u_m^{(2)}&v_{m}^{(2)}
\end{pmatrix}
.
\end{eqnarray}
Note that $\varphi_n$ and $\psi_n$ satisfy the equation $-\frac{d^2}{dx^2}f+Vf=0$ which is equivalent to 
\[
\begin{pmatrix}
\frac{d}{dx}f\\
\frac{d^2}{dx^2}f
\end{pmatrix}
=
\begin{pmatrix}
0&1\\
V&0
\end{pmatrix}
\begin{pmatrix}
f\\
\frac{d}{dx}f
\end{pmatrix}.
\]
Differentiating both sides of $(\ref{derivative of constants})$, we obtain 
\begin{eqnarray*}
\begin{pmatrix}
0&0\\
V&0
\end{pmatrix}
\begin{pmatrix}
\tilde{\varphi}_{m}&\tilde{\psi}_{m}\\
\frac{d}{dx}\tilde{\varphi}_{m}&\frac{d}{dx}\tilde{\psi}_{m}
\end{pmatrix}
\begin{pmatrix}
u_{m}^{(1)}&v_{m}^{(1)}\\
u_m^{(2)}&v_{m}^{(2)}
\end{pmatrix}
=
\begin{pmatrix}
\tilde{\varphi}_{m}&\tilde{\psi}_{m}\\
\frac{d}{dx}\tilde{\varphi}_{m}&\frac{d}{dx}\tilde{\psi}_{m}
\end{pmatrix}
\begin{pmatrix}
\frac{d}{dx}u_{m}^{(1)}&\frac{d}{dx}v_{m}^{(1)}\\
\frac{d}{dx}u_m^{(2)}&\frac{d}{dx}v_{m}^{(2)}
\end{pmatrix}.
\end{eqnarray*}
Thus we have 
\begin{eqnarray}\label{deribative of constants2}
\begin{pmatrix}
\frac{d}{dx}u_{m}^{(1)}&\frac{d}{dx}v_{m}^{(1)}\\
\frac{d}{dx}u_m^{(2)}&\frac{d}{dx}v_{m}^{(2)}
\end{pmatrix}
=
-V
\begin{pmatrix}
\tilde{\varphi}_{m}\tilde{\psi}_{m}&\tilde{\psi}_{m}^2\\
-\tilde{\varphi}_{m}^2&-\tilde{\varphi}_{m}\tilde{\psi}_{m}
\end{pmatrix}
\begin{pmatrix}
u_{m}^{(1)}&v_{m}^{(1)}\\
u_m^{(2)}&v_{m}^{(2)}
\end{pmatrix}.
\end{eqnarray}
Let $u_m=
\begin{pmatrix}
u_m^{(1)}\\
u_m^{(2)}
\end{pmatrix}$ and $v_m=
\begin{pmatrix}
v_m^{(1)}\\
v_m^{(2)}
\end{pmatrix}$. 
By $(\ref{deribative of constants2})$, we see that $u_m$ and $v_m$ satisfy that 
\begin{eqnarray}\label{derivative of constants3}
\frac{d}{dx}u_m&=&
-V
\begin{pmatrix}
\tilde{\varphi}_{m}\tilde{\psi}_{m}&\tilde{\psi}_{m}^2\\
-\tilde{\varphi}_{m}^2&-\tilde{\varphi}_{m}\tilde{\psi}_{m}
\end{pmatrix}
u_m,
\\
\label{derivative of constants4}
\frac{d}{dx}v_m&=&
-V
\begin{pmatrix}
\tilde{\varphi}_{m}\tilde{\psi}_{m}&\tilde{\psi}_{m}^2\\
-\tilde{\varphi}_{m}^2&-\tilde{\varphi}_{m}\tilde{\psi}_{m}
\end{pmatrix}
v_m.
\end{eqnarray}
\begin{Lemma}\label{wronskian constancy}
For $x\in I_m$, 
$u_{m}^{(1)}(x)v_{m}^{(2)}(x)-v_{m}^{(1)}(x)u_m^{(2)}(x)=1$.
\end{Lemma}
\begin{Proof}\rm
By $(\ref{derivative of constants3})$ and $(\ref{derivative of constants4})$, we obtain
\begin{eqnarray}\nonumber
\frac{d}{dx}(u_{m}^{(1)}v_{m}^{(2)}-v_{m}^{(1)}u_m^{(2)})&=&
\frac{d}{dx}u_{m}^{(1)}v_{m}^{(2)}+u_{m}^{(1)}\frac{d}{dx}v_{m}^{(2)}-\frac{d}{dx}v_{m}^{(1)}u_m^{(2)}-v_{m}^{(1)}\frac{d}{dx}u_m^{(2)}
\\
&=&\nonumber
-V(\tilde{\varphi}_{m}\tilde{\psi}_{m}u_m^{(1)}+\tilde{\psi}_{m}^2u_m^{(2)})v_{m}^{(2)}
+u_{m}^{(1)}V(\tilde{\varphi}_{m}^2v_m^{(1)}+\tilde{\varphi}_{m}\tilde{\psi}_{m}v_m^{(2)})
\\
&&\nonumber
+V(\tilde{\varphi}_{m}\tilde{\psi}_{m}v_m^{(1)}+\tilde{\psi}_{m}^2v_m^{(2)})u_m^{(2)}
-v_{m}^{(1)}V(\tilde{\varphi}_{m}^2u_m^{(1)}+\tilde{\varphi}_{m}\tilde{\psi}_{m}u_m^{(2)})
\\
&=&0.\nonumber
\end{eqnarray}
Since $u_{m}^{(1)}(x_m-\alpha_m)v_{m}^{(2)}(x_m-\alpha_m)-v_{m}^{(1)}(x_m-\alpha_m)u_m^{(2)}(x_m-\alpha_m)=1$, we have our assertion.
\qed
\end{Proof}
\begin{Lemma}\label{lower bound of constants matrix}
Assume $u_j,v_j \in \mathbb{R}$, $j=1,2$ and 
$u_1v_2-v_1u_2=1$. 
Then 
\[
\low
\begin{pmatrix}
u_1&v_1\\
u_2&v_2
\end{pmatrix}
\geq \cfrac{1}{\sqrt{u_1^2+u_2^2+v_1^2+v_2^2}}.
\]
\end{Lemma}
\begin{Proof}\rm
Let $u=
\begin{pmatrix}
u_1\\
u_2
\end{pmatrix}$ and  
$v=
\begin{pmatrix}
v_1\\
v_2
\end{pmatrix}
$, and 
$u\cdot v=u_1v_1+u_2v_2$.
The eigenvalues $\xi_{\pm}$ of the matrix $
\begin{pmatrix}
u_1&u_2\\
v_1&v_2
\end{pmatrix}
\begin{pmatrix}
u_1&v_1\\
u_2&v_2
\end{pmatrix}
$ are 
\[
\xi_{\pm}=\frac{1}{2}
\left(
|u|^2+|v|^2\pm \sqrt{
(|u|^2+|v|^2)^2-4(|u|^2|v|^2-u\cdot v)^2
}
\right).
\]
Since $1-\sqrt{1-t}\geq\frac{1}{2}t$ for $0<t<1$ and $|u|^2|v|^2-u\cdot v=(u_1v_2-u_2v_1)^2=1$, we obtain 
\[
\xi_{-}\geq\frac{1}{|u|^2+|v|^2}.
\]
This implies our assertion.
\qed
\end{Proof}
%We estimate $\low R_m.$
\begin{Lemma}
It follows that 
\begin{eqnarray}\label{lower bound of R_m}
\low R_m \geq 
\cfrac{1}{2\sqrt{2\alpha_m^2+1}}
\exp
\left(
-\frac{1}{3}h_m(4\alpha_m^3+3\alpha_m)
\right).
\end{eqnarray}
\end{Lemma}
\begin{Proof}\rm
Note that
\[
R_m=\lim_{x\uparrow x_m+\alpha_m}
\begin{pmatrix}
\varphi_{m}(x)&\psi_{m}(x)\\
\frac{d}{dx}\varphi_{m}(x)&\frac{d}{dx}\psi_{m}(x)
\end{pmatrix}=
\lim_{x\uparrow x_m+\alpha_m}
\begin{pmatrix}
\tilde{\varphi}_{m}(x)&\tilde{\psi}_{m}(x)\\
\frac{d}{dx}\tilde{\varphi}_{m}(x)&\frac{d}{dx}\tilde{\psi}_{m}(x)
\end{pmatrix}
\begin{pmatrix}
u_{m}^{(1)}(x)&v_{m}^{(1)}(x)\\
u_m^{(2)}(x)&v_{m}^{(2)}(x)
\end{pmatrix}.
\]
It is straightfoward to see
\begin{eqnarray}
\lim_{x\uparrow x_m+\alpha_m}
\begin{pmatrix}
\tilde{\varphi}_{m}(x)&\tilde{\psi}_{m}(x)\\
\frac{d}{dx}\tilde{\varphi}_{m}(x)&\frac{d}{dx}\tilde{\psi}_{m}(x)
\end{pmatrix}
&=&\nonumber
\begin{pmatrix}
1&2\alpha_m\\
0&1
\end{pmatrix},
\\
\low
\begin{pmatrix}
1&2\alpha_m\\
0&1
\end{pmatrix}
&\geq& \cfrac{1}{\sqrt{4\alpha_m^2+2}}.
\label{lower bound of non-pertubed matrix}
\end{eqnarray}
%\left(\begin{array}{cc}u_{m}^{(1)}(x)&v_{m}^{(1)}(x)\\u_m^{(2)}(x)&v_{m}^{(2)}(x)\end{array}\right). 
By Lemmas $\ref{wronskian constancy}$ and $\ref{lower bound of constants matrix}$, we have 
\begin{eqnarray*}
\low 
\begin{pmatrix}
u_{m}^{(1)}(x)&v_{m}^{(1)}(x)\\
u_m^{(2)}(x)&v_{m}^{(2)}(x)
\end{pmatrix}
\geq
\frac{1}{\sqrt{|u_m(x)|^2+|v_m(x)|^2}}.
\end{eqnarray*}
We see that for $a, b \in \mathbb{R}$, 
\begin{eqnarray}
\sup_{\theta\in[0,\pi)}
\left|
\left(
\begin{pmatrix}
\cos \theta\\
\sin \theta
\end{pmatrix}
,
\begin{pmatrix}
ab&b^2\\
-a^2&-ab
\end{pmatrix}
\begin{pmatrix}
\cos \theta\\
\sin \theta
\end{pmatrix}
\right)
\right|=\frac{a^2+b^2}{2}.
\label{q-norm}
\end{eqnarray}
By $(\ref{derivative of constants3})$ and $(\ref{q-norm})$, we have 
\begin{eqnarray}\nonumber
\frac{d}{dx}(|u_m|^2)&=&2\left(u_m, \frac{d}{dx} u_m\right)\\
&\leq&\nonumber
2|V|\left|
\left(u_m,
\begin{pmatrix}
\tilde{\varphi}_{m}\tilde{\psi}_{m}&\tilde{\psi}_{m}^2\\
-\tilde{\varphi}_{m}^2&-\tilde{\varphi}_{m}\tilde{\psi}_{m}
\end{pmatrix}u_m \right)\right|
\\ \nonumber
&\leq&
h_m(\tilde{\varphi}_{m}^2+\tilde{\psi}_{m}^2)|u_m|^2.
\end{eqnarray}
Thus, by Gronwall's inequality, we obtain
\begin{eqnarray*}
|u_m(x)|^2
\leq
\exp
\left(
h_m
\int_{x_m-\alpha_m}^{x}
\left(\tilde{\varphi}_{m}(y)^2+\tilde{\psi}_{m}(y)^2\right)
dy
\right).
\end{eqnarray*}
In particular
\begin{eqnarray}\nonumber
|u_m(x_m+\alpha_m)|^2
&\leq&
\exp
\left(
h_m
\int_{x_m-\alpha_m}^{x_m+\alpha_m}
\left(1+(y-x_m+\alpha_m)^2\right)
dy
\right).
\\
\nonumber
&=&
\exp
\left(
h_m
(2\alpha_m+\frac{8}{3}\alpha_m^3)
\right).
\end{eqnarray}
We can estimate $|v_m|^2$ in a similar way:
\begin{eqnarray*}
|v_m(x_m+\alpha_m)|^2
\leq
\exp
\left(
\frac{2}{3}h_m(4\alpha_m^3+3\alpha_m)
\right).
\end{eqnarray*}
Therefore, we see that 
\begin{eqnarray}
\low 
\left.
\begin{pmatrix}
u_{m}^{(1)}(x)&v_{m}^{(1)}(x)\\
u_m^{(2)}(x)&v_{m}^{(2)}(x)
\end{pmatrix}
\right|_{x=x_m+\alpha_m}
&\geq&
\left.
\frac{1}{\sqrt{|u_m(x)|^2+|v_m(x)|^2}}
\right|_{x=x_m+\alpha_m}
\nonumber
\\
\label{lower bound of constants}
&\geq&
\cfrac{1}{\sqrt{2}}
\exp
\left(
-\frac{1}{3}h_m(4\alpha_m^3+3\alpha_m)
\right).
\end{eqnarray}
By $(\ref{lower bound of non-pertubed matrix})$ and $(\ref{lower bound of constants})$, we have our assertion.
\qed
\end{Proof}
\begin{proof_mainT}\rm
By $(\ref{lower bound of W_m})$, $(\ref{lower bound of R_m})$, and $c_n=R_nW_{n-1}c_{n-1}$ for $n\geq1$,
we obtain
\begin{eqnarray}\label{lower bound of c_n}
|c_n|\geq 
\frac{1}{2^n}
\left(
\prod_{m=1}^{n}(L_{m-1}^2+2)(2\alpha_m^2+1)
\right)^{-\frac{1}{2}}
\exp
\left(
-\frac{1}{3}\sum_{m=1}^{n}h_m(4\alpha_m^3+3\alpha_m)
\right).
\end{eqnarray}
If $n$ is sufficiently large, then 
\begin{eqnarray*}
\frac{L_n^4}{L_n^3+3L_n}\geq
\frac{L_n}{4}.
\end{eqnarray*}
Thus, by $(\ref{lower bound of L^2 integral})$ and $(\ref{lower bound of c_n})$, if $n$ is sufficiently large, then we have
\begin{eqnarray*}
\int_{J_n}|f(x)^2|dx
\geq
\frac{L_n}{4^{n+2}}
\left(
\prod_{m=1}^{n}(L_{m-1}^2+2)(2\alpha_m^2+1)
\right)^{-1}
\exp
\left(-\frac{2}{3}\sum_{m=1}^{n}h_m(4\alpha_m^3+3\alpha_m)
\right).
\end{eqnarray*} 
Suppose $(\ref{assumption})$. Then we obtain 
\begin{eqnarray}\nonumber
\int_{J_n}|f(x)^2|dx\rightarrow \infty.
\end{eqnarray}
This implies that any solutions $f$ of $-\frac{d^2}{dx^2}f+Vf=0$ do not belong to $L^2([0,\infty))$.
\qed
\end{proof_mainT}
\section{Proof of Corollary 1.3.}
Note that $x_n=\exp(n^n)$, $\alpha_n=\frac{1}{2}$, $h_n=e^{n}$, $L_0=x_1-\frac{1}{2}$ and $L_n={x_{n+1}}-{x_n}-1$ for $n\geq1$.

\begin{Lemma}\label{AAA}
If there exists $f \in D(H_{\theta})$ such that $(f,H_{\theta}f)<0$, then $H_{\theta}$ has a single negative eigenvalue.
\end{Lemma}
\begin{Proof}\rm
$\mathcal{R}[A]$ denotes the range of a map $A$. 
We see that $\dim\mathcal{R}[E_{\theta}((-\infty,0))]\leq1$. 
For its proof, see Lemma \ref{appendix}. 
Let $f \in D(H_{\theta})$ such that $(f,H_{\theta}f)<0$. 
Then this implies $\mathcal{R}[E_{\theta}((-\infty,0))] \neq \{0\}$ and $\dim\mathcal{R}[E_{\theta}((-\infty,0))]=1$.
This implies our assertion.
\qed
\end{Proof}
\begin{Lemma}\label{BBB} 
It follows that 
\begin{enumerate}[$(1)$]
\item
$H_{\theta}$ has no negative eigenvalues for $\theta \in [0,\frac{\pi}{2}]$,
\item
$H_{\theta}$ has a single negative eigenvalue for $\theta \in (-\frac{\pi}{2},\arctan(-\frac{1+\sqrt{3}}{2})]$.
\end{enumerate}
\end{Lemma}
\begin{Proof}\rm
Let $f \in D(H_{\theta})$. 
Then $f$ satisfies the boundary condition $f(0)\sin\theta-\frac{d}{dx}f(0)\cos\theta=0$.
We obtain that 
\begin{eqnarray*}
(f,H_{\theta}f)&=&\int_{0}^{\infty}\overline{f(x)}(-\frac{d^2}{dx^2}f(x)+V(x)f(x))dx\\
&=&
\overline{f(0)}\frac{d}{dx}f(0)+
\int_{0}^{\infty}\left(\left|\frac{d}{dx}f(x)\right|^2+V(x)|f(x)|^2\right)dx.
\end{eqnarray*}
If $\theta=0$ or $\frac{\pi}{2}$, then we have $\overline{f(0)}\frac{d}{dx}f(0)=0$ and $(f,H_{\theta}f)\geq0$. 
If $0<\theta<\frac{\pi}{2}$, then by the boundary condition, we obtain
\begin{eqnarray*}
(f,H_{\theta}f)=
|f(0)|^2\tan\theta+\int_{0}^{\infty}\left(\left|\frac{d}{dx}f(x)\right|^2+V(x)|f(x)|^2\right)dx
\geq0.
\end{eqnarray*}
Thus $(f,H_{\theta}f)\geq0$ for $\theta\in[0,\frac{\pi}{2}]$. 
This implies the first part of our assertion. 

We shall prove that there exists $f \in D(H_{\theta})$ such that $(f,H_{\theta}f)<0$ for any $\theta \in (-\frac{\pi}{2},\arctan(-\frac{1+\sqrt{3}}{2})]$. 
It is sufficient to prove there exists $f\in L^2([0,\infty))$ such that $(-\frac{d^2}{dx^2}f+Vf,f)<0$ for any boundray conditions $\frac{\frac{d}{dx}f(0)}{f(0)}=-\lambda$, $\lambda\geq\frac{1+\sqrt{3}}{2}$. 
Let $\lambda\geq1$. 
Define $f_{\lambda}:[0,\infty)\rightarrow \mathbb{R}$ by 
\[
f_{\lambda}(x)=
\begin{cases}
\exp\left( \cfrac{\lambda}{x-1} \right),& \text{if $0\leq x <1$}\\
0,&\text{otherwise}.
\end{cases}
\]
We see that $f_{\lambda} \in L^2([0,\infty))$ and that $f_{\lambda}(0)=\exp(-\lambda)$, $\frac{d}{dx}f_{\lambda}(0)=-\lambda\exp(-\lambda)$. 
%Let $\delta_{\lambda}=\sup\{t^2\exp(-\lambda t) \mid t\geq1 \}.$
Then we have 
\begin{eqnarray*}
\int_{0}^{\infty}|\frac{d}{dx}f_{\lambda}(x)|^2dx&=&
\int_0^{1}
\left|
\cfrac{\lambda}{(x-1)^2}\exp \left(\frac{\lambda}{x-1} \right) 
\right|^2dx
\\
&=&
\frac{1}{4\lambda}(2\lambda^2+2\lambda+1)\exp(-2\lambda).
\end{eqnarray*}
%We see that 
%\[
%\delta_{\lambda}=
%\begin{cases}
%\exp(-\lambda),& \text{if $\lambda\geq2$},\\
%4\exp(-2) \lambda^{-2},& \text{if $0 \leq \lambda < 2$}.
%\end{cases}
%\]
Since $V(x)f_{\lambda}(x)=0$ for $x\geq0$, we obtain that, for 
$\lambda \geq \frac{1+\sqrt{3}}{2}$,
\begin{eqnarray*}
\left(-\frac{d^2}{dx^2}f_{\lambda}+Vf_{\lambda},f_{\lambda}\right)
=\frac{1}{4\lambda}(-2\lambda^2+2\lambda+1)\exp(-2\lambda)
<0.
\end{eqnarray*}
By Lemma $\ref{AAA}$, we have our assertion.
\qed
\end{Proof}
\begin{Lemma}\label{lem}
For any $p>0$, it follows that
\begin{eqnarray*}
\lim_{n\rightarrow \infty}
x_n
\left(
\prod_{m=1}^{n-1}x_m
\right)^{-p}=\infty.
\end{eqnarray*}
\end{Lemma}
\begin{Proof}\rm
We obtain 
\begin{eqnarray*}
x_{n+1}
\left(
\prod_{m=1}^{n}x_m
\right)^{-p}
&=&
\exp
\left(
(n+1)^{n+1}-pn^n-p\sum_{m=1}^{n-1}m^m
\right)
\\
&\geq&
\exp
\left(
(n+1)^{n+1}-pn^n-p(n-1)^n
\right)
\\
&=&
\exp
\left(
(n+1-2p)(n+1)^{n}
\right)
\\
&\rightarrow&\infty,\qquad\text{as $n\rightarrow \infty$.}
\end{eqnarray*}
\qed
\end{Proof}
\begin{proof_main}\rm
\noindent
By Lemma $\ref{BBB}$ it is sufficient to prove that $0\notin \sigma_{pp}(H_\theta)$ for all $\theta$.
We see that for all $n\geq 1$,
\begin{eqnarray*}
\frac{L_n^2+2}{L_n^2}<2 .
\end{eqnarray*}
Thus we have for all $n\geq 1$,
\begin{eqnarray*}
\left(
\prod_{m=1}^{n}(L_{m-1}^2+2)
\right)^{-1}
\geq
\frac{2^{n-1}}{L_0^2+2}
\left(
\prod_{m=1}^{n-1}L_m
\right)^{-2}.
\end{eqnarray*}
We see that $L_n<x_{n+1}$ for all $n\geq1$ and $L_n>\frac{1}{2}x_{n+1}$ for sufficiently large $n\geq1$. 
Therefore we have 
\begin{eqnarray}
L_n
\left(
\prod_{m=1}^{n-1}L_m
\right)^{-2}
\geq
\cfrac{1}{2}
x_{n+1}
\left(
\prod_{m=1}^{n}x_m
\right)^{-2}.
\nonumber
\end{eqnarray} 
By Lemma $\ref{lem}$, we obtain
\begin{eqnarray}
&&
\cfrac{L_n}{4^n}
\left(
\prod_{m=1}^{n}(L_{m-1}^2+2)
\right)^{-1}
\left(
\prod_{m=1}^{n}(2\alpha_m^2+1)
\right)^{-1}
\exp
\left(-\frac{2}{3}\sum_{m=1}^{n}h_m(4\alpha_m^3+3\alpha_m)
\right)
\nonumber
\\
&&\geq
\frac{1}{2(L_0^2+2)3^{n}}L_n
\left(
\prod_{m=1}^{n-1}L_m
\right)^{-2}
\exp
\left(
-\frac{4}{3}
\sum_{m=1}^{n}e^m
\right)
\nonumber
\\
&&\geq
\frac{1}{4(L_0^2+2)3^{n}}x_{n+1}
\left(
\prod_{m=1}^{n}x_m
\right)^{-2}
\exp
\left(
-\frac{4}{3}e^{n+1}
\right)
\nonumber
\\
&&=
\frac{\exp(\frac{1}{3}(n+1)^{n+1})}{4(L_0^2+2)3^{n}}
{x_{n+1}}^{\frac{1}{3}}
\left(
\prod_{m=1}^{n}x_m
\right)^{-2}
\exp
\left(\frac{1}{3}
\left(
(n+1)^{n+1}
-4e^{n+1}
\right)
\right)
\nonumber
\\
&&
\rightarrow \infty,\qquad\text{as $n\rightarrow \infty$.}
\nonumber
\end{eqnarray}
By Theorem $\ref{main}$, we see that $0\notin \sigma(H_\theta)$.
\qed
\end{proof_main}

\begin{appendices}
\section{}
By \cite{Simonsparse}, we see that $\sigma_{ac}(H_{\theta})=\emptyset$ and $\sigma_{pp}(H_{\theta})\cap(0,\infty)=\emptyset$. 
In this appendix, we prove that $H_{\theta}$ has a single negative eigenvalue for some $\theta$ and that $\sigma_{sc}(H_{\theta})=[0,\infty)$ for all $\theta$. 
Let $V$ ba a sparse potential with 
$x_n=\exp(n^n)$ for $n=1,2,...$, and 
\begin{eqnarray*}
V(x)=
\begin{cases}
e^{n},&\text{if $|x-x_n|\leq\frac{1}{2}$ for $n=1,2,...$},\\
0,&\text{otherwise}.
\end{cases}
\end{eqnarray*}
%Let $H=-\frac{d^2}{dx^2}+V:L^2([0,\infty))\rightarrow L^2([0,\infty))$ with its domain $C^{\infty}_0((0,\infty))$ and $H_{\theta}$ be an self-adjoint extension of $H$. 
By \cite{Simonsparse} and \cite{SimonSpencer}, we see that $\sigma_{pp}(H_{\theta})\cap(0,\infty)=\emptyset$ and $\sigma_{ac}(H_{\theta})=\emptyset$ for all $\theta$. 
%The next lemma says that $H_{\theta}$ has a single negative eigenvalue depending on $\theta$.
\begin{Lemma}\label{appendix}
Let $E_{\theta}$ be the spectral measure of $H_{\theta}$. 
For all $\theta$, $\dim \mathcal{R}[ E_{\theta}((-\infty,0))]\leq1$.
\end{Lemma}
\begin{Proof}\rm
We prove this by a contradiction. 
Suppose that $\dim \mathcal{R}[ E_{\theta}((-\infty,0))]>1$. 
Then we can take $\varphi,\psi \in E_{\theta}((-\infty,0))$ such that $\varphi$ and $\psi$ are orthogonal to each other. 
Let $\varphi_n=E_{\theta}((-n,-\frac{1}{n}))\varphi$ and $\psi_n=E_{\theta}((-n,-\frac{1}{n}))\psi$. 
We see that $\varphi_n,\psi_n\in D(H_{\theta})$, $\varphi_n \rightarrow \varphi$ and $\psi_n\rightarrow \psi$. 
Let $N\geq1$ be sufficiently large such that $\varphi_N$ and $\psi_N$ are linearly independent. Since $\alpha\varphi_N+\beta\psi_N\in \mathcal{R}[E_{\theta}((-N,-\frac{1}{N}))]$ for $\alpha,\beta\in \mathbb{C}$, we have, for $(\alpha,\beta) \neq (0,0)$, 
\begin{eqnarray}\label{condition1}
(\alpha\varphi_N+\beta\psi_N,H_{\theta}(\alpha\varphi_N+\beta\psi_N))<0.
\end{eqnarray}

On the other hand, since the deficiency indices of $H$ are equal to one, there exists an ismometric operator $U_{\theta}:\ker(H^*-i)\rightarrow \ker(H^*+i)$ and $w \in \ker(H^*-i)$ such that 
\[
D(H_\theta)=\left\{
v+\alpha(w+U_{\theta}w)\middle| v\in D(\overline{H}), \alpha \in \mathbb{C} 
\right\},
\]
where $\overline{H}$ is the closure of $H$.
Let $u_\theta=w+U_\theta w$.
There exist $v_1,v_2 \in D(\overline{H})$ and $\alpha_1,\alpha_2\in \mathbb{C}$ such that 
\begin{eqnarray*}
\varphi_N&=&v_1+\alpha_1u_{\theta},\\
\psi_N&=&v_2+\alpha_2u_{\theta}.
\end{eqnarray*}
If $\alpha_1=0$, then $(\varphi_N,H_{\theta}\varphi_N)=(v_1,\overline{H}v_1)\geq0$. 
$(\ref{condition1})$ implies $\alpha_1\neq0$. 
Similarly we have $\alpha_2 \neq 0$.
We obtain 
\begin{eqnarray*}\label{condition2}
(\alpha_2\varphi_N+\alpha_1\psi_N,H_{\theta}(\alpha_2\varphi_N+\alpha_1\psi_N))
=(\alpha_2v_1+\alpha_1v_2,\overline{H}(\alpha_2v_1+\alpha_1v_2))\geq0.
\end{eqnarray*}
This contradicts with $(\ref{condition1})$.
Thus we have our assertion.
\qed
\end{Proof}Let $E<0$. 
We see that $\dim \ker(H^*-E)=1$. 
This implies that there exists a boundary condition $\theta(E)$ such that $H_{\theta(E)}$ has a single negative eigenvalue $E$.
%The next lemma gives singular continous spectrum of $H_{\theta}$. 
\begin{Lemma}
For all $\theta \in (-\frac{\pi}{2}, \frac{\pi}{2}]$, it follows that $\sigma_{sc}(H_{\theta})=[0,\infty)$. 
\end{Lemma}
\begin{Proof}\rm
Since $\sigma_{pp}(H_{\theta})\cap(0,\infty)=\emptyset$ and $\sigma_{ac}(H_{\theta})=\emptyset$ for all $\theta$, 
we see that $\sigma(H_{\theta})\cap(0,\infty)=\sigma_{sc}(H_{\theta})\cap(0,\infty)$. 
We prove $(0,\infty) \subset \sigma(H_{\theta})$ for all $\theta$ by contradition. 
Suppose that there exist $\theta\in[0,\pi),E>0$ such that $E \in (0,\infty)\setminus \sigma(H_{\theta})$.
Since $H$ is regular at zero and limit-circle case at infinity, the deficiency indices $\dim \ker(H^*\pm i)$ are equal to one. 
Thus $\dim \ker (H^*-E)=1$. 
This implies that there exists an $L^2$-solution of $-\frac{d^2}{dx^2}f+Vf=Ef$. 
By \cite[Theorem 2.3.]{Simonsparse}, however, $-\frac{d^2}{dx^2}f+Vf=pf$ has no solutions with $f \in L^2([0,\infty))$ for any $p>0$. 
This is a contradiction and we get $(0,\infty)\subset\sigma_{sc}(H_{\theta})$ for all $\theta$.
Since $\sigma_{sc}(H_{\theta})\cap(-\infty,0)=\emptyset$, we get our assertion.
\qed
\end{Proof}
\section*{Acknowledgement}
This work was supported by JST SPRING, Grant Number JPMJSP2136.

\end{appendices}

\section*{Declarations}
\subsection*{Ethical Approval }
The author has no competing interests to declare that are relevant to the content of this article.

\subsection*{Competing interests}
The author declares no conflict of interest.

\subsection*{Authors' contributions }
The author confirms sole responsibility for this manuscript .

\subsection*{Funding}
This work was supported by JST SPRING, Grant Number JPMJSP2136.

\subsection*{Availability of data and materials }
Data sharing not applicable to this article as no datasets were generated or analysed during the current study.

\end{document}